\documentclass{amsproc}

\usepackage{color, hyperref, amssymb}
\usepackage{mathrsfs}

\newtheorem{theorem}{Theorem}[section]
\newtheorem{lemma}[theorem]{Lemma}

\theoremstyle{definition}
\newtheorem{definition}[theorem]{Definition}
\newtheorem{example}[theorem]{Example}

\newtheorem{problem}{PROBLEM}

\theoremstyle{remark}

\numberwithin{equation}{section}

\newcommand{\fin}{\text{\rm fin}}
\newcommand{\red}{{\text{\rm red}}}
\newcommand{\N}{\mathbb N}
\newcommand{\Z}{\mathbb Z}

\newcommand{\DP}{\negthinspace :\negthinspace}

\DeclareMathOperator{\spec}{spec} 
\DeclareMathOperator{\Pic}{Pic} 
 
\DeclareMathOperator{\Kr}{Kr}

\begin{document}

\title[On Long-Term Problems]{On Long-Term Problems in \\ Multiplicative Ideal Theory and Factorization Theory}

\author{Alfred Geroldinger}
\address{Department of Mathematics and Scientific Computing, University of Graz, \newline 8010 Graz, Austria}
\email{alfred.geroldinger@uni-graz.at}
\thanks{The first author was supported  by the Austrian Science Fund FWF, Project P36852-N}

\author{Hwankoo Kim}
\address{Division of Computer Engineering, Hoseo University, Asan 31499, \newline Republic of Korea}
\email{hkkim@hoseo.edu}
\thanks{The second author was supported in part by the Basic Science Research Program through the National Research Foundation of Korea (NRF), funded by the Ministry of Education (2021R1I1A3047469).}

\author{K.~Alan Loper}
\address{Department of Mathematics, The Ohio State University, Columbus. OH 43210, USA}
\email{lopera@math.ohio-state.edu}




\begin{abstract}
In this survey article we discuss key open problems which could serve as a guidance for further research directions of multiplicative ideal theory and factorization theory.
\end{abstract}

\maketitle

\section{Introduction} \label{sec:introduction}
\smallskip

The origins of multiplicative ideal theory can be traced back to the beginnings of commutative algebra in the first half of the 20th century. We refer to articles of Halter-Koch \cite{HK98, HK25a} for a discussion of the early development. The monographs of Gilmer \cite{Gi72} and Larsen-McCarthy \cite{{La-Mc71}}, which appeared in the 1970s,  strongly influenced  the area. Building on the work of Lorenzen and Aubert, Halter-Koch showed (in his book \cite{HK98}, published 1998) that large parts of multiplicative ideal theory can be derived in the purely multiplicative setting of commutative monoids. This viewpoint also turned out to be crucial for the further development of factorization theory.

Factorization theory has its origin in algebraic number theory. Starting points were  results  (in the 1960s) of Carlitz and of Narkiewicz, who studied the asymptotic behavior of arithmetic counting functions in orders of algebraic number fields. From the 1990s on, factorization theory expanded its scope into commutative ring theory and since then, multiplicative ideal theory and factorization theory developed in parallel and stimulated each other.
At present, both areas  are in the midst of a stormy development. This is documented by
recent research monographs of Chabert, Elliott, Fontana et al., Geroldinger et al., Grynkiewicz, Halter-Koch, Knebusch et al., Wang and Kim (\cite{Kn-Zh02a, Kn-Ka14a, Fo-Ho-Lu13, Ch25a, El19a,Gr22a, HK25a,Wa-Ki24a,Ge-Gr-Zh25a}), by a steady flow of conferences devoted to these fields and the associated conference proceedings (for some of the last 10 years, see \cite{FFG14, CFGO16, FFGTZ17, FFGO20, BC21, CFFGJ23, B-G-O-S25}).

 \newpage
The above mentioned monographs contain an abundance of open problems (in addition, we  would like to refer to the survey \cite{C-F-F-G14} by Cahen et al. on unsolved questions). Thus we made no attempt to give a comprehensive overview on open problems, but we tried to offer brief discussions of some research directions which we hope will stimulate further research in multiplicative ideal and factorization theory.

\section{Preliminaries} \label{sec:preliminaries}

Let $\N$ denote the set of positive integers and set $\N_0 = \N \cup \{0\}$. For integers $a, b \in \Z$, let $[a, b] = \{ x \in \Z \colon a \le x \le b \}$ denote the discrete interval between $a$ and $b$. For subsets $A, B$ of an additive abelian group $G$, let $A+B = \{a+b \colon a \in A, b \in B\} \subseteq G$ denote their sumset.
If  $S$ is any set, then $\mathcal P (S)$ denotes the system of all subsets of $S$ and $\mathcal P_{\fin} (S)$ the system of all finite, nonempty subsets of $S$.

All semigroups are commutative and have an identity element. If not stated otherwise, we use multiplicative notation. Let $S$ be a semigroup. We denote by $S^{\times}$ the group of invertible elements, and we say that $S$ is {\it reduced} if $S^{\times} = \{1_S\}$. Then $S_{\red} = S/S^{\times} = \{aS^{\times} \colon a \in S\}$ denotes the associated reduced semigroup of $S$. We say that $S$ is {\it cancellative} if $a, b, c \in S$ and $ab=ac$ implies that $b=c$.
By a {\it monoid}, we mean a  cancellative semigroup. Let $M$ be a monoid. An element $a \in M$ is called an {\it atom} (or {\it irreducible}) if $a \notin M^{\times}$ and an equation $a=bc$, with $b, c \in M$, implies that $b \in M^{\times}$ or $c \in M^{\times}$. We denote  by $\mathcal A (M)$ the set of atoms  and by $\mathsf q (M)$ the quotient group of $M$.

\smallskip
Next we introduce a key concept of multiplicative ideal theory, namely ideal systems or, equivalently,  star-operations ($*$-operations). The term star-operation was popularized by Gilmer, following the term {\it Strich-operation} introduced by Krull. Following Aubert, Halter-Koch used the term {\it ideal system}. These concepts were generalized to semistar operations and then to abstract closure operations by Elliott. Ideal systems were generalized to weak ideal systems and then to module systems by Halter-Koch \cite{HK11b, HK25a}. In the present survey we use
both the terms star-operation and ideal system.

\smallskip
Let $M$ be a monoid.
An {\it ideal system} on $M$ is a map $r\colon\mathcal P(M)\to\mathcal P(M)$, defined by $X \mapsto X_r := r (X)$ for all $X \subseteq M$,  such that the following conditions are satisfied for all subsets $X,Y\subseteq M$ and all $c\in M$:
\begin{itemize}
\item $X\subseteq X_r$.

\item $X\subseteq Y_r$ implies $X_r\subseteq Y_r$.

\item $cM\subseteq\{c\}_r$.

\item $cX_r=(cX)_r$.
\end{itemize}
Let $r$ be an ideal system on $M$. A subset $I \subseteq M$ is called an {\it $r$-ideal} if $I_r=I$. We denote by $\mathcal I_r (M)$ the set of all nonempty $r$-ideals, and we define $r$-multiplication by setting $I \cdot_r J = (IJ)_r$ for all $I,J \in \mathcal I_r (M)$. Then $\mathcal I_r (M)$ together with $r$-multiplication is a reduced semigroup with identity element $M$, and it is a multiplicative lattice. Conversely, by \cite{Du-Ep-Gi25}, for every multiplicative lattice there are a commutative semigroup $H$ and a weak ideal system $r$ on $H$ such that $L$ is isomorphic to the lattice of $r$-ideals of $H$.

A subset $X \subseteq \mathsf q (M)$ is called ($M$-){\it fractional} if $cX \subseteq M$ for some $c \in M$. We define
\[
X_r := c^{-1} (cX)_r \,,
\]
and note that the definition does not depend on the choice of $c$. A subset $I \subseteq \mathsf q (M)$ is a {\it fractional $r$-ideal} if $I = X_r$ for some fractional subset $X \subseteq \mathsf q (M)$. We denote by
 $\mathcal F_r (M)$ the set of fractional $r$-ideals. For $I, J \in \mathcal F_r (M)$, the set $IJ$ is $M$-fractional, and we call $I \cdot_r J := (IJ)_r \in \mathcal F_r (M)$ the $r$-product of $I$ and $J$. A fractional $r$-ideal $J \in \mathcal F_r (M)$ is called \emph{$r$-invertible} if there exists a fractional $r$-ideal $J'$ such that $J \cdot_r J' = M$.  $(\mathcal F_r (M), \cdot_r)$ is a  semigroup containing $\mathcal I_r (M)$, $\mathcal F_r (M)^{\times}$ is the group of $r$-invertible fractional $r$-ideals, and $\mathcal I_r^*(M) = \mathcal F_r^{\times} (M) \cap \mathcal I_r (M)$ the  monoid of $r$-invertible $r$-ideals of $M$ with $r$-multiplication.

The map $\delta \colon \mathsf q (M) \to \mathcal F_r (M)^{\times}$, defined by $\delta (a) = aM$, is a group homomorphism with kernel $\ker ( \delta ) = M^{\times}$. Its cokernel
\[
\mathcal C_r (M) = \mathcal F_r (M)^{\times}/ \delta ( \mathsf q (M))
\]
is called the {\it $r$-class group} of $M$. For subsets $X, Y \subseteq \mathsf q(M)$, we denote by $X_s = XM$, by $(X \DP Y)=\{x\in\mathsf q(M)\colon xY \subseteq X\}$, by $X^{-1}=(M \DP X)$, by $X_v=(X^{-1})^{-1}$ and by $X_t=\bigcup_{E\subseteq X,|E|<\infty} E_v$.
We will use the $s$-system, the $v$-system, the $t$-system, and the $w$-system. The $v$-system was introduced by Krull and $v$-ideals are also called divisorial ideals. The $w$-system was introduced by R.~McCasland and F.~Wang (see Section \ref{sec:w-ideals}).

\smallskip
All rings are commutative and have an identity element.
For a commutative ring $D$, we denote by $D^{\times}$ its group of units and by
$D^{\bullet}$ its monoid of regular elements. Let $D$ be a domain with quotient field $K$. We denote by $\mathcal A (D)$ the set of atoms, by $\mathsf q (D)$ the quotient field of $D$,  by $\overline D$ the integral closure of $D$, and by $\widehat D = \widehat{D^{\bullet}} \cup \{0\}$ the complete integral closure of $D$.

An {\it ideal system} on $D$ is a map $r\colon\mathcal P(D)\to\mathcal P(D)$ such that the following conditions are satisfied for all subsets $X,Y\subseteq D$ and all $c\in D$:
\begin{itemize}
\item $X \cup \{0\} \subseteq X_r$.

\item $X\subseteq Y_r$ implies $X_r\subseteq Y_r$.

\item $cM\subseteq\{c\}_r$.

\item $cX_r=(cX)_r$.
\end{itemize}
If $r \colon\mathcal P(D)\to\mathcal P(D)$ is an ideal system on $D$, then $r^{\bullet} \colon\mathcal P(D^{\bullet})\to\mathcal P(D^{\bullet})$, defined by $r^{\bullet} (X) = r (X) \setminus  \{0\}$ for every $X \subseteq D$, is an ideal system on $D^{\bullet}$. Conversely, if $r \colon\mathcal P(D^{\bullet})\to\mathcal P(D^{\bullet})$ is an ideal system on $D^{\bullet}$, then $r^{\circ} \colon \mathcal P(D)\to\mathcal P(D)$, defined by $r^{\circ} (X) = r (X) \cup \{0\}$ for every $X \subseteq D$, is an ideal system on $D$. Thus, in the sequel we will not distinguish between $r, r^{\bullet}$, and $r^{\circ}$ and simply speak of $r$-ideals of $D$ and of $D^{\bullet}$. In particular, $\mathcal I_r (D)$, $\mathcal F_r (D)$, $\mathcal I^*_r (D)$, and $\mathcal C_r (D) = \mathcal C_r (D^{\bullet})$ have their obvious meaning.

The $d$-system $d \colon \mathcal P(D)\to\mathcal P(D)$, defined by $X_d = \langle X \rangle$ for every $X\subseteq D$, is the system of classical ring ideals of $D$. Whenever we speak of ideals of $D$, we mean $d$-ideals and in all notation introduced above we omit the $d$-suffix. In particular, $\mathcal I^* (D) = \mathcal I_d^* (D)$ is the monoid of invertible ideals of $D$ and $\mathcal C_d (D) = \Pic (D)$ is the Picard group of $D$.

A monoid (resp. a domain) is said to be
\begin{itemize}
\item {\it Mori} (resp. {\it $v$-Noetherian}) if it satisfies the ascending chain condition on $v$-ideals.

\item {\it strong Mori}  (resp. {\it $w$-Noetherian})  if it satisfies the ascending chain condition on $w$-ideals.

\item\label{Krull-charact-1} {\it Krull} if it is completely integrally closed and Mori.
\end{itemize}
Every $v$-ideal is a $t$-ideal and if the  monoid is Mori, then $v=t$. Krull monoids and Noetherian domains are strong Mori, and
every strong Mori monoid is Mori, because every $t$-ideal is a $w$-ideal.

\section{Intersecting Valuation Domains } \label{sec:intersection}

Let $D$ be a domain with quotient field $K$ and suppose that $D \subsetneq K$.
In general, there are a lot of domains between $D$
and $K$.  For example, if $D$ is a Noetherian domain with Krull
dimension greater than one, then there are uncountably many domains between
$D$ and $K$.

A valuation domain $V$ with quotient field $L$ is, in a sense, close
to being all of $L$.  For example, every domain between $V$ and $L$
is again a valuation domain and is a localization of $V$.  And if $V$ has finite Krull dimension
then there are only finitely many such intermediate domains.

A theorem of Krull should make the study of  valuation domains worthwhile.

\begin{theorem}
Every integrally closed domain $D$ (not a field) can be expressed as
the intersection of all the valuation domains between $D$ and its quotient field.
\end{theorem}

It seems reasonable then that valuation domains should play a central
role in the study of integrally closed domains.  They are the
indivisible building blocks from which everything is made.  Some
work has been done in this direction but not a lot.  A principal
difficulty is that given an integrally closed domain $D$
it is usually hard to determine what the valuation domains are.
And, given a collection of valuation domains it is hard to discern
what the structure of the intersection is.  A valuation overring
$V$ of $D$ is local, and its maximal ideal contracts to a prime
ideal of $D$.  But, often, the maximal ideals of uncountably many
valuation overrings all contract to the same prime ideal of $D$.

Clearly, many efforts are made to build new examples of domains as
overrings of well-understood domains.  This is generally done by
starting with a domain that is well-understood and adding in some
fractions to obtain an overring.  Instead of this bottom-up approach we explore here a top-down approach of intersecting valuation domains.

\begin{problem} \label{intersect-valuation-1}~

\begin{enumerate}
\item Given an integrally closed domain what can we say about the
collection of its valuation overrings?

\item Given a collection of valuation domains what can we say about
the domain that results from intersecting them?
\end{enumerate}
\end{problem}

These are both difficult questions about which little is known.
In this note we focus on one facet of the second one.
As we observed above, valuation domains are local.  One could then
hope that, given a domain, we could obtain valuation overrings
by localizing at prime ideals.  This rarely succeeds.  But sometimes
it does.  We say that a domain is a Pr\"ufer domain if the localization
at any prime ideal is a valuation domain.  So, if given a Pr\"ufer
domain, it is easy to write the domain as an intersection of
valuation domains.  The representation is simply the intersection
of the localizations at prime ideals.  The other direction is
a fascinating question.

 \begin{problem}  \label{intersect-valuation-2}
 Given a collection of valuation domains, what are necessary and sufficient conditions
that we can put on the collection that would make the intersection a
Pr\"ufer domain?
\end{problem}

The following theorem gives partial answers to Problem \ref{intersect-valuation-2}.

\begin{theorem}
Let  $\{V_i \ | \ i \in \Omega \} $ be a  collection of
valuation domains that all have the same quotient field.  In each of the following cases the
intersection of the valuation domains is a Pr\"ufer domain.
\begin{enumerate}
\item $\Omega$ is finite.  \cite{M-Na62}

\item There is a  monic non-constant
polynomial  $f \in K[X]$, whose coefficients lie in each $V_i$ and
which has no root in the residue field of any $V_i$.  \cite{Gi69}

\item There is a field $L$, whose cardinality  is strictly larger than the cardinality
of $\Omega$,  such that $L \subseteq V_i$ for each $i \in \Omega$.  \cite{OR06}
\end{enumerate}
\end{theorem}

The common theme of these three properties (1) - (3) seems to be that the collection
of valuation domains needs to be, in some sense, thin or sparse.  This
observation, in turn, might lead one to seek a topological condition.
We digress to describe a construction of Krull which shows that looking
at the standard topologies on a collection of valuation domains is
not sufficient to determine when the intersection is a Pr\"ufer
domain.

We start by giving a quick demonstration that Pr\"ufer domains have
properties that make them easy to work with.

Let $K$ be a field and let $D = K[X,Y]$ be the
polynomial ring in $X$ and $Y$.  Let $J = (X,Y)$ and
$I = (X^2, Y^2)$ be ideals of $D$.  It is easy to see that
$J^3 = J I$.  It seems intuitively clear then that
$J^2 = I$ should be true.  But it is not.  On the other hand,
suppose that $T$ is a Pr\"ufer domain and let $d,r \in T$ be nonzero.
Let $B = (d,r)$ and $H = (d^2, r^2)$ be ideals of $T$.  As above,
$BH = B^3$. is clear.  Now, however, in the context of a Pr\"ufer
domain, we can cancel $B$ and conclude that $B^2 = H$ is true.

Krull's idea was to define certain closure operations on the collection of fractional ideals of a domain.  These closure operations are known
as star operations.  Of particular interest are the star operations
he identified as ab (arithmetisch brauchbar = arithmetically
useful).  Let $D$ be an integrally closed domain and let $*$ be an
ab star operation.  Then when $I,J,H$ are finitely generated ideals
such that $(IH)^* \subseteq (JH)^*$ then we can conclude that $I^* \subseteq J^*$.

Before pressing on, we recap what we have discussed so far and where
we want to go with it.

\begin{itemize}
\item Every integrally closed domain is the intersection of its
valuation overrings.

\item Pr\"ufer domains are precisely the domains for which any
localization at a prime ideal is a valuation domain.

\item We would like to know what conditions on a collection of valuation
rings with a common quotient field characterize those collections
for which the intersection is a Pr\"ufer domain.

\item We have some results which seem to indicate that the property
(or properties) that work will be topological in nature.
\end{itemize}

There are two well-known topologies which can be placed on a collection
of valuation domains (for a sample of contributions in this direction  see \cite{Co-Ol26a, Du-Go23a,Fi-Go-Sp23,Fi-Lo24a,Go24a}).  These are known as
\begin{enumerate}
\item the Zariski topology, and

\item the constructible or patch or ultrafilter topology.
\end{enumerate}
Some amazing work of Krull demonstrates clearly that neither of these
topologies alone can provide a solution to our problem.

Let $V$ be a valuation domain.   For a polynomial $f \in V[X]$, let $c(f)$  be the ideal
of $V$ generated by the coefficients of $f$.
Let $N = \{ f \in V[X] \mid c(f) = V\}$ and finally, let
$V(X) = N^{-1} V[X]$.  $V(X)$ is known as the {\it Nagata extension} of $V$.
$V(X)$ is a valuation domain which is almost identical to $V$.  The only
significant difference between $V$ and $V(X)$ is that $V(X)$ has more
units.

Let $D$ be an integrally closed domain (not a field) with quotient field $K$.  Krull showed
how to use ab star operations on $D$ to build what he called
Kronecker function domains.  \cite{FL06}  There is always a unique minimal
Kronecker function domain which we label as $\Kr(D)$ which has
the following properties.

\begin{enumerate}
\item $\Kr(D)$ is an overring of $D[X]$.

\item $\Kr(D) \cap K = D$.

\item $\Kr(D)$ is a Pr\"ufer domain.

\item The localizations $\Kr(D)_{\mathfrak p}$ where $\mathfrak p$ runs over the prime ideals
of $\Kr(D)$ are precisely the valuation domains $V(X)$, where $V$ runs
over the valuation overrings of $D$.
\end{enumerate}

The key importance of $\Kr(D)$ for us is that when $D$ is an integrally
closed domain the space of valuation overrings of $D$ and the space of valuation overrings of $\Kr(D)$ are homeomorphic for both of the topologies noted above.  This is problematic because, while $D$ can be any integrally closed domain, $\Kr(D)$ is always a Pr\"ufer domain.  Hence, there is no hope of characterizing when an intersection of valuation domains is a Pr\"ufer domain using just information about the
topology of the collection of valuation domains.

\section{Classification of Integral Domains between $\mathbb{Z}[X]$ and $\mathbb{Q}[X]$}\label{kim_sec_kim_prob1}

We survey recent developments in the classification of  domains that lie strictly between the polynomial ring $\mathbb{Z}[X]$ and $\mathbb{Q}[X]$. We focus on two major studies:
\begin{itemize}
    \item[\cite{L-T09}] which classifies \textit{integrally closed domains} (including Pr\"ufer and Dedekind domains) via order-theoretic arguments, extending classical valuation-theoretic methods.
    \item[\cite{Pe25a}] which classifies \textit{Dedekind domains} using the concept of \textit{generalized rings of integer-valued polynomials}, broadening earlier results on factorization and class groups.
\end{itemize}

Both works build on classical ideas from valuation theory, factorization theory, and the arithmetic of polynomial rings, offering new insights into how certain local constructions (valuation overrings) glue together to form global integrally closed domains between $\mathbb{Z}[X]$ and $\mathbb{Q}[X]$.

The primary motivation of \cite{L-T09} is to give an order-theoretic classification of integrally closed domains situated between $\mathbb{Z}_{(p)}[X]$ and $\mathbb{Q}[X]$, where $p$ is a prime number and $\mathbb{Z}_{(p)}$ denotes the localization of $\mathbb{Z}$ at $p\mathbb{Z}$. In particular, \cite{L-T09} focuses on the structure of valuation overrings and extends classical results of MacLane on constructing valuations via \textit{key polynomials}.

A representative theorem from \cite{L-T09} is:

\begin{theorem}[{\cite{L-T09}}]
\label{LT09_theorem}
Let $p \in \mathbb{Z}$ be a fixed prime, and let $D$ be a domain such that
\[
\mathbb{Z}_{(p)}[X] \;\subseteq\; D \;\subseteq\; \mathbb{Q}[X].
\]
Then $D$ is Pr\"ufer if and only if $D$ is integrally closed and all the $p$-unitary valuation overrings of $D$ are limit.
\end{theorem}

Constructive aspects hinge on the notion of \textit{key polynomials}, which allow for systematic building of valuation overrings, ultimately classifying the possible integrally closed domains in a step-by-step manner.

Having explored the order-theoretic classification of integrally closed domains, we now shift our focus to Dedekind domains. To facilitate this discussion, we first introduce the relevant notation.

Let $\mathbb{P}$ be the set of all prime numbers. For a fixed $p \in \mathbb{P}$, we adopt the following notation:

\begin{itemize}
    \item $\mathbb{Z}_p$ denotes the ring of $p$-adic integers.
    \item $\mathbb{Q}_p$ denotes the field of $p$-adic integers.
    \item $\mathbb{C}_p$ denotes the completion of the algebraic closure of the field $\mathbb{Q}_p$.
    \item $\overline{\mathbb{Z}}_p$ denotes the absolute integral closure of $\mathbb{Z}_p$.
    \item $\overline{\widehat{\mathbb{Z}}} = \prod_{p \in \mathbb{P}} \overline{\mathbb{Z}}_p$.
\end{itemize}

A \textit{generalized ring of integer-valued polynomials} is defined, for a subset $\underline{E} \subset \overline{\widehat{\mathbb{Z}}}$, as
\[
\mathrm{Int}_\mathbb{Q}(\underline{E},\overline{\widehat{\mathbb{Z}}}) \;=\; \{\,f\in \mathbb{Q}[X] \,\mid\, f(\alpha)\in \overline{\widehat{\mathbb{Z}}} \;\text{for all}\; \alpha\in \underline{E}\}.
\]
This idea generalizes the classical ring of integer-valued polynomials.

Let $\underline{E} = \prod_p E_p$ be a subset of $\overline{\widehat{\mathbb{Z}}}$. Following \cite[Definition 2.11]{Pe23a}, we say that $\underline{E}$ is \textit{polynomially factorizable} if, for every $g \in \mathbb{Z}[X]$
and every $\alpha = (\alpha_p) \in \underline{E}$, there exist $n, d \in \mathbb{Z}$ with $n, d \geq 1$ such that $g(\alpha)^n / d$ is a unit of $\overline{\widehat{\mathbb{Z}}}$.

In \cite{Pe23a}, the author classifies \textit{Dedekind domains with finite residue fields of prime characteristic} that lie between $\mathbb{Z}[X]$ and $\mathbb{Q}[X]$. These domains are precisely realized as generalized rings of integer-valued polynomials. The work is inspired by earlier studies on factorization theory, integer-valued polynomials, and class group computations.

A key observation is that if $W$ is a one-dimensional valuation overring of such a Dedekind domain with $\mathbb{Q} \subset W$, then $W = \mathbb{Q}[X]_{(q)}$ for some irreducible polynomial $q \in \mathbb{Q}[X]$, implying that its residue field has infinite cardinality. In \cite{Pe25a}, the author extends this classification by relaxing the condition on the cardinality of residue fields of prime characteristic.

Recall that an element $\alpha \in \mathbb{C}_p$ has \emph{bounded ramification} if the extension $\mathbb{Q}_p \subseteq \mathbb{Q}_p(\alpha)$ has finite ramification. Denote by $\mathbb{C}_p^{br}$ the set of all elements of $\mathbb{C}_p$ with bounded ramification, and let $\mathbb{O}_p$ be the unique local ring of $\mathbb{C}_p$. Define $\mathcal{O} = \prod_{p \in \mathbb{P}} \mathbb{O}_p \subset \prod_{p \in \mathbb{P}} \mathbb{C}_p$.

\begin{theorem}[{\cite[Theorem 3.4]{Pe25a}}]
\label{Pe25a_theorem}
Let $R$ be a Dedekind domain such that $\mathbb{Z}[X] \subset R \subseteq \mathbb{Q}[X]$. Then $R$ is equal to $\operatorname{Int}_{\mathbb{Q}}(\underline{E}, \mathcal{O})$ for some polynomially factorizable subset $\underline{E} = \prod_{p \in \mathbb{P}} E_p \subset \mathcal{O}$, where for each prime $p$, the set $E_p \subset \mathbb{O}_p \cap \mathbb{C}_p^{br}$ consists of finitely many transcendental elements over $\mathbb{Q}$.
\end{theorem}

This result completes the classification of Dedekind domains between $\mathbb{Z}[X]$ and $\mathbb{Q}[X]$ and provides a powerful characterization of Dedekind domains in this setting, linking local valuation data and integer-valued polynomial constraints to global structural properties such as the class group. The result is obtained by characterizing discrete valuation rings (DVRs) $W$ of the field $\mathbb{Q}(X)$ that are residually algebraic over $\mathbb{Z}_{(p)}$ for some prime $p$ (see \cite[Corollary 2.28]{Pe25a}).

When the extension of residue fields has infinite degree, such a domain is given by
\[
\{\phi \in \mathbb{Q}(X) \mid v_p(\phi(\alpha)) \geq 0\}
\]
for some $\alpha \in \mathbb{C}_p$, which is transcendental over $\mathbb{Q}_p$. The case where the residue field extension has finite degree was addressed in \cite{Pe17a}.

This result also generalizes findings by Loper, Werner, and others concerning Prüfer domains and valuation-theoretic constructions.

Although \cite{L-T09} and \cite{Pe23a, Pe25a} offer significant advances in classifying integrally closed domains between $\mathbb{Z}[X]$ and $\mathbb{Q}[X]$, the broader landscape of such domains is still only partially explored. In particular, unifying \textit{generalized integer-valued polynomial} viewpoints with \textit{valuation-theoretic} frameworks remains an active area of research, as do possible applications to:
\begin{itemize}
\item \emph{Krull or $($strong$)$ Mori domains}: Extending classification results to settings with additional constraints (e.g., finite residue fields, specific $v$-class groups, etc.).

\item \emph{Factorial domains}: Investigating when generalized integer-valued polynomial rings are factorial.
\end{itemize}

Recall that an integral domain $D$ is called an \emph{almost Dedekind domain} (resp., a \emph{$t$-almost Dedekind domain}) if $D_{\frak m}$ is a DVR for all maximal ideals (resp., maximal $t$-ideals) $\frak{m}$ of $D$.

\begin{problem}\label{kim_prob1}
Classify all integral domains of a specified type (e.g., ($t$-)almost Dedekind domains and Krull domains satisfying additional conditions) and determine their placement between $\mathbb{Z}[X]$ and $\mathbb{Q}[X]$. Furthermore, investigate the extent to which such domains can be realized as:
\begin{itemize}
    \item Generalized rings of integer-valued polynomials, building on \cite{Pe23a, Pe25a}, or
    \item Valuation overrings constructed via MacLane-type key polynomial techniques, as in \cite{L-T09}.
\end{itemize}
\end{problem}

\section{Generalizing Kronecker Function Domains} \label{sec:kronecker}

In a previous section we discussed a construction of Krull known as the Kronecker function domain.  We start with an integrally closed
domain $D$ with quotient field $K$ and build the domain $\Kr(D)$ which satisfies

\begin{enumerate}

\item $\Kr(D)$ lies between $D[X]$ and $K(X)$.

\item $\Kr(D)$ is a Pr\"ufer domain.

\item The prime ideals of $\Kr(D)$ correspond precisely with the valuation overrings of $D$.

\end{enumerate}

There are two features of this that might strike one as annoyances that one would wish to avoid.  We start as above with an integrally
closed domain $D$ with quotient field $K$.

\begin{itemize}
\item The actual construction is built using the ab star operations.  It is cumbersome.  One could wish for a simpler
construction.

\item One could wonder whether the base ring $D$ is essential.  Is it possible to build a Kronecker function domain using
just $K$ - without $D$?
\end{itemize}

In 2003 Halter-Koch   \cite{HK03}     gave what he called an axiomatic treatment of what he called Kronecker  domains which gives us the construction
we were looking for above.

Let $K$ be a field which has valuation subrings.  Let $X$ be an indeterminate.  We say that a ring $T$ is a {\it Kronecker domain} provided

\begin{enumerate}
\item $T \subseteq K(X)$.

\item $X$ is a unit in $T$.

\item For every nonzero polynomial $f \in K[X]$, we have $f(0) \in fT$.
\end{enumerate}

If these conditions are satisfied, then we have also

\begin{itemize}
\item $T$ is a Pr\"ufer domain.

\item If $\mathfrak p$ is a prime ideal of $T$ then $T_{\mathfrak p}$ has the form $V(X)$ for some
valuation domain $V$ contained in $K$.

\item The conventional Kronecker function domains defined for subrings $D$ of $K$
are included as Kronecker domains.
\end{itemize}

Halter-Koch's notion of Kronecker domain has been known since 2003 and very little
has been done in that time to study them.   Actually, Halter-Koch himself extended what is known by himself by demonstrating that they can be built with a variant of the conventional construction (\cite[Theorem 6.5.5]{HK25a}).
Many questions could be posed.  We will content ourselves with one.

\begin{problem}
Fabbri and Heubo-Kwegna \cite{Fa-HK12a} give explicit details for the construction of a non-traditional Kronecker domain in a very
specialized setting.  Can their results be generalized to a much wider setting?
\end{problem}

\section{Coherent Domains} \label{sec:coherent}

Coherent domains were first defined by  Chase \cite{Ch60} in 1960.
 A domain $D$ is {\it coherent} if it satisfies one of the  following three equivalent
 conditions (note there are many more equivalent ones).

 \begin{enumerate}
\item[(a)] The intersection of any two finitely generated ideals of $D$
is finitely generated.

\item[(b)] Any finitely generated ideal of $D$ is finitely presented.

\item[(c)] The direct product of any family of  flat $D$-modules is flat.
 \end{enumerate}

We focus on the first of these conditions.  It is clear from this
condition that any Noetherian domain is coherent.  It is not as
simple, but still straightforward to prove that any Pr\"ufer
domain is coherent.  When one is considering all domains there are
other natural examples of coherent domains.  In the realm of one-dimensional domains things are much more restrictive.
A one-dimensional coherent domain which is integrally closed
is a Pr\"ufer domain.
The integral closure of a one-dimensional Noetherian domain is
a Dedekind domain which is both Noetherian and Pr\"ufer.
The class of one-dimensional coherent domains which are neither
Noetherian nor Pr\"ufer is vanishingly small.  We are aware of
one such example in the literature.  This example is not
integrally closed but has Pr\"ufer integral closure.

This state of affairs led Vasconcelos to ask

\begin{problem}
Is the integral closure of a one-dimensional coherent domain a Pr\"ufer domain? \cite{G92}
\end{problem}

One might hope to make progress on this question by considering
many examples of one-dimensional coherent domains that are neither
integrally closed nor Noetherian.  But again, such examples
scarcely exist.

Instead of looking at this question from a forward direction
(examine the integral closure of a one-dimensional coherent
domain) maybe one can make progress by looking at the question
from a backwards perspective.

\begin{itemize}
\item Let $D$ be a one-dimensional integrally closed domain
    which is not Pr\"ufer.  Is $D$ in any sense close to being
    coherent?
\end{itemize}

A recent result shows that if such a domain $D$ is local it is,
in a sense anti-coherent.

\begin{theorem}
Let $D$ be a one-dimensional integrally closed local domain which is not a Pr\"ufer
domain.  Let $I,J$ be two principal ideals of $D$ with no containment
relations.  Then $I \cap J$ is not finitely generated.  \cite{GL21}
\end{theorem}

This seems to be evidence toward a positive answer to Vasconcelos'
question.  How could a coherent domain - where the intersection of
two finitely generated ideals is always finitely generated - have as
its integral closure
a domain where the intersection of two principal ideals is
almost never finitely generated?

\begin{problem}
Let $K$ be a field.  Let $X,Y$ be indeterminates.  Consider the
valuation domain $V = K(X)[Y]_{(Y)}$.  Let $M$ be the maximal ideal of
$V$.  Let $D = K + M$.  Then $D$ is a one-dimensional local domain
which is integrally closed and not Pr\"ufer.  Can we prove that
$D$ is not the integral closure of a coherent domain?
\end{problem}

\section{Theory of $w$-modules} \label{sec:w-ideals}

Originally formulated by McCasland and Wang \cite{Wa-Mc97}, \cite{M-W99}, \emph{$w$-module theory} refines classical module and ideal theories in commutative algebra by focusing on a specific star operation, the $w$-operation. The concept of $w$-ideals
coincides with the $F_\infty$-operation previously considered by J.R.~Hedstrom and E.G.~Houston \cite{He-Ho80}, and the $w$-ideals coincide with the semidivisorial ideals of S.~Glaz and W.~Vasconcelos \cite{Gl-Va77}).

The framework of $w$-ideals generalizes earlier results on $t$-ideals and divisorial ideals while providing new insights in non-Noetherian settings. The core motivation for $w$-module theory arose from studying the structure of commutative rings and the interplay of their star operations with modules and ideals.

\begin{definition}\label{def:GVwtheory}~
Let $R$ be a commutative ring, possibly with zero divisors.
\begin{enumerate}
\item A nonzero ideal $J$ of $R$ is called a \emph{Glaz--Vasconcelos ideal} (GV-ideal) if $J$ is finitely generated and the natural homomorphism
\[
\varphi \colon R \;\longrightarrow\; \operatorname{Hom}_{R}(J, R)
\]
is an isomorphism. The set $\operatorname{GV}(R)$ of all GV-ideals forms a multiplicative system of ideals in $R$.

\item An $R$-module $M$ is called \emph{GV-torsion-free} if, for any $J \in \operatorname{GV}(R)$, the condition $Jx = 0$ (for some $x \in M$) implies $x = 0$.
On the other hand, $M$ is called \emph{GV-torsion} if, for every nonzero element $x \in M$, there exists $J \in \operatorname{GV}(R)$ such that $Jx = 0$.
It is evident that the pair (GV-torsion modules, GV-torsion-free modules) defines a hereditary torsion theory, first introduced in \cite[p. 1618]{ki-li-Zh19}.

\item For a GV-torsion-free $R$-module $M$, define
\[
M_{w} \;=\; \bigl\{\,x \in E(M) \;\bigm|\; \text{there exists }J \in \operatorname{GV}(R)\text{ such that }Jx \subseteq M\bigr\},
\]
where $E(M)$ is the injective envelope of $M$. The module $M_{w}$ is called the \emph{$w$-envelope} of $M$.

\item A GV-torsion-free $R$-module $M$ is called a \emph{$w$-module} if $\operatorname{Ext}_{R}^{1}(R/J, M) = 0$ for every $J \in \operatorname{GV}(R)$; equivalently, if $M = M_w$. In particular, an ideal $I$ of $R$ is called a \emph{$w$-ideal} if $I = I_w$.

\item An $R$-module $M$ is said to be \emph{$w$-flat} if, for every maximal $w$-ideal $\mathfrak{m}$ of $R$, the localization $M_{\mathfrak{m}}$ is flat over $R_{\mathfrak{m}}$.
\end{enumerate}
\end{definition}

\noindent
The $w$-operation builds a bridge towards homological algebra. The resulting $w$-module theory can be used to study certain ``high-dimensional" domains from a module-theoretic perspective. In fact, the $w$-module theory is a type of locally homological method.

Every flat $R$-module is $w$-flat, and if $D$ is a Pr\"ufer domain, $w$-flat $D$-modules coincide with flat $D$-modules. This shows that $w$-module theory specializes to classical notions in more restricted settings.

\smallskip
It is well-known that Noetherian rings play an important role in commutative algebra and algebraic geometry. Using $w$-operations, strong Mori domains were introduced in \cite{Wa-Mc97}.
It was previously claimed (\cite[Corollary~IX.5.5]{Fu-Sa01}) that for a domain $D$ with quotient field $K$, the injective module $E(K/D)$ is $\Sigma$-injective if and only if $D$ is a Mori domain. However, this statement is \emph{incorrect}: one must replace `Mori' with `strong Mori' (see \cite[Theorem~9]{Fu03}).
In recent years, many classical theorems about Noetherian rings have been generalized to strong Mori domains.

In parallel with $w$-module theory, the study of \emph{Krull domains} has been pivotal in commutative algebra via the lens of star operations. Indeed, a Krull domain $D$ can be characterized by the following conditions.
\begin{enumerate}
\item[(i)] For every height-one prime ideal $\mathfrak{p}$, the localization $D_{\mathfrak{p}}$ is a discrete valuation domain.
\item[(ii)] $D = \bigcap D_{\mathfrak{p}}$, the intersection taken over all height-one prime ideals $\mathfrak{p}$.
\item[(iii)] Every nonzero non-unit in $D$ is contained in only finitely many height-one prime ideals.
\end{enumerate}
Analogous to how Dedekind domains can be characterized by the invertibility of nonzero ideals, Krull domains can be characterized as follows (\cite[Theorem~3.6]{Ka89a}, \cite[Theorem~7.9.3]{Wa-Ki24a}):
\[
  D \text{ is Krull}
  \quad \Longleftrightarrow \quad
  \text{every nonzero ideal of $D$ is $t$- (or $w$-) invertible}.
\]

This implies that a Krull domain can essentially be considered a $t$- (or $w$-) Dedekind domain. Furthermore, it suggests that Krull domains can be studied in a manner similar to Dedekind domains, utilizing star-operator methods.

\begin{problem}\label{kim_prob2}
Prove or disprove the assertion that
\begin{itemize}
\item[] \emph{Every result in classical commutative algebra has a corresponding analogue in $w$-module theory}.
\end{itemize}
\end{problem}

Finding an accurate analogue for this may be challenging and could require the development of new concepts and techniques. For instance, in a Dedekind domain, any finitely generated module $M$ decomposes as a direct sum of a projective module and a torsion module. A central question in $w$-module theory (and more broadly in star-operation settings) is whether such classical decompositions have analogues in Krull domains. For Krull domains, what is the corresponding formulation for finitely generated ($w$-finitely generated) modules?

Recent extensions explore $w$-factor rings. For a domain $D$, the following conditions are equivalent.
\begin{enumerate}
    \item[(a)] $D$ is Krull.
    \item[(b)] The $w$-factor ring $(D/I)_w$ is a zero-dimensional semilocal principal ideal ring for every proper $w$-ideal $I$ of $D$, \cite[Theorem~4.5]{Ch-Ki23a}.
    \item[(c)] $(D/I)_w$ is a quasi-Frobenius ring \cite[Theorem 9]{Xi-Zh24}.
\end{enumerate}
 These results highlight a strong parallel between Dedekind and Krull domains when analyzed through $w$-operations and $w$-ideal structures. This raises the question of whether P$v$MDs can be characterized by their $w$-factor rings, where an integral domain $D$ is called a \emph{Pr\"ufer $v$-multiplication domain} (\emph{P$v$MD} for short) if every nonzero finitely generated ideal of $D$ is $t$-invertible (or $w$-invertible).

\begin{problem}
Investigate whether the P$v$MDs $D$ are precisely those domains for which their $w$-factor rings $(D/I)_w$ are FC-rings for every proper $w$-ideal $I$ of $D$, where a commutative ring $R$ is an \emph{FC-ring} if it is coherent and self-FP-injective.
\end{problem}

\section{Homological versus Elementary Methods} \label{sec:homological}

\smallskip
\centerline{\it Every ring considered in this section is an integral domain.}
\smallskip

In 1937 Krull defined a regular local ring to be a Noetherian local
ring $D$ with maximal ideal $\mathfrak m$ for which the minimum number of elements in a generating set for $\mathfrak m$ is the same as the height of $\mathfrak m$  \cite{Kr37}.
This definition is completely elementary.

In 1956 Serre  \cite{Se56}    proved that
a Noetherian local ring is regular if and only if it has finite global
dimension.  This homological characterization of regular local rings
was pursued because it led to an easy proof that a localization of a
regular local ring at a prime ideal was still regular, which had
proven to be quite elusive using non-homological methods.

This event is often cited as validating the use of homological methods
in the study of ring theory.  In fact, homological methods have
become quite dominant.  We propose several settings in which the
elementary methods might regain prominence.

\smallskip
{\it  Question I.}  Does there now exist a non-homological proof that
any localization of a regular local ring is still regular?

Answer:  There exist informal sites on the internet that indicate
that the answer is yes.  These sites indicate that  the proof involves careful use of regular
sequences. This is not surprising since it is well-known that
an ideal generated by a regular sequence has finite projective
dimension.  But we were unable to find a published proof.  So

\begin{problem}  Can one prove non homologically that any
localization of a regular local ring is still regular?
\end{problem}

We push this a little farther.  Ratliff  \cite{Ra71} defined
a set $\{d_0, d_1, \ldots , d_n \}$ of elements in a local domain
$(D, \mathfrak m)$ to be analytically independent if, whenever
$F(d_0, \ldots, d_n) = 0$ for some homogeneous polynomial
$F(X_0, \ldots, X_n) \in D[X_0, \ldots, X_n]$, then all coefficients of
$F$ are in $\mathfrak m$.  This generalizes the notion of a regular
sequence.  Say that the analytically independent dimension
is the maximal size of an analytically independent set.  Is it possible that the notion of finite projective
dimension could be replaced by the notion of analytically independent
dimension?

\begin{problem}
An ideal generated by a regular sequence has finite projective dimension.  A regular sequence is
analytically independent.  Is it true that an ideal generated by an analytically independent set has
finite projective dimension?
\end{problem}

\smallskip
{\it Question II}:  Can we use homological terminology to define
classes of rings that are not being studied at all and might be
profitably studied with elementary techniques?

Answer: Perhaps.  Recall that a Dedekind domain is a domain for
which every nonzero ideal is invertible.  The more general notion
of a Pr\"ufer domain is a domain for which every finitely
generated ideal is invertible.  We could similarly say that a domain
has pseudo-finite global dimension provided there is an upper
bound on the projective dimensions of all finitely generated modules.  Perhaps the union of an ascending
tower of regular local rings could provide nontrivial examples
of such rings.  Examples of such towers are given in
recent papers concerning the so-called Shannon extensions  \cite{He17}.

\begin{problem}
Can the property of having pseudo-finite global dimension be defined in elementary terms?
\end{problem}

\begin{problem}
Do Shannon extensions have pseudo-finite global dimension?
\end{problem}

\section{Gorenstein Multiplicative Ideal Theory} \label{sec:gorenstein}

Gorenstein homological algebra arose as an extension of classical homological algebra to better accommodate \emph{non-regular} contexts. Traditional projective, injective, and flat modules often fail to capture subtler behaviors in more general rings and modules, especially those beyond Noetherian or regular situations. Gorenstein dimensions (and related Gorenstein notions) were introduced to address these shortcomings, building on Auslander and Bridger's foundational work on Gorenstein projective modules (\cite{Au-Br69}). Their investigations were initially guided by syzygy theory over Noetherian rings and the concept of complete intersections.

Subsequent expansions by Enochs, Jenda, and Torrecillas introduced a range of new tools---such as relative homological dimensions (\cite{En-Je95, En-Je-To93}). Classical homological algebra heavily emphasizes projective, injective, and flat modules. The Gorenstein versions of these notions `relax' certain regularity assumptions.

\smallskip
\centerline{\it Throughout this section,}
\centerline{\it   let $R$ denote a commutative ring, possibly with zero divisors.}
\smallskip

\begin{definition}
Let $M$ be an $R$-module.
\begin{enumerate}
    \item $M$ is said to be \emph{Gorenstein projective} if there exists an exact sequence of projective $R$-modules
    \[
    \cdots \; \longrightarrow \; P_1 \; \longrightarrow \; P_0 \; \longrightarrow \; P^0 \; \longrightarrow \; P^1 \; \longrightarrow \; \cdots
    \]
    such that $M \cong \ker(P^0 \to P^1)$ and $\operatorname{Hom}_R(-, Q)$ preserves exactness of this sequence for every projective $R$-module $Q$.

    \item $M$ is said to be \emph{Gorenstein injective} if there exists an exact sequence of injective $R$-modules
    \[
    \cdots \; \longrightarrow \; I_1 \; \longrightarrow \; I_0 \; \longrightarrow \; I^0 \; \longrightarrow \; I^1 \; \longrightarrow \; \cdots
    \]
    such that $M \cong \ker(I^0 \to I^1)$ and $\operatorname{Hom}_R(Q, -)$ preserves exactness of this sequence for every injective $R$-module $Q$.

    \item $M$ is said to be \emph{Gorenstein flat} if there exists an exact sequence of flat $R$-modules
    \[
    \cdots \; \longrightarrow \; F_1 \; \longrightarrow \; F_0 \; \longrightarrow \; F^0 \; \longrightarrow \; F^1 \; \longrightarrow \; \cdots
    \]
    such that $M \cong \ker(F^0 \to F^1)$ and exactness is preserved under the functor $-\otimes_R E$ for every injective $R$-module $E$.

    \item The \emph{Gorenstein projective dimension} of $M$, denoted $\operatorname{Gpd}_R(M)$, is defined as the smallest integer $n$ (if such an $n$ exists) such that there exists an exact sequence of Gorenstein projective modules:
    \[
        0 \longrightarrow G_n \longrightarrow G_{n-1} \longrightarrow \cdots \longrightarrow G_1 \longrightarrow G_0 \longrightarrow M \longrightarrow 0,
    \]
    where each $G_i$ is a Gorenstein projective module. If no such finite $n$ exists, we set $\operatorname{Gpd}_R(M) = \infty$.
\end{enumerate}
\end{definition}

\medskip

\begin{example}~

\begin{enumerate}
\item If $R$ is a \emph{Gorenstein ring}, then every finitely generated $R$-module has finite Gorenstein projective dimension. Over a regular local ring, the Gorenstein projective dimension agrees with the usual projective dimension.

\item Let $R = K[X]/(X^2)$ for a field $K$. This is a Frobenius algebra, so $R$ is both Gorenstein projective and Gorenstein injective as an $R$-module.
\end{enumerate}
\end{example}

\begin{definition}
An integral domain $D$ is called a
\begin{itemize}
\item \emph{Gorenstein Dedekind domain} (briefly, G-Dedekind domain) if \emph{every ideal of $D$ is G-projective}.

\item \emph{Gorenstein Pr\"ufer domain} (briefly, G-Pr\"ufer domain) if $D$ is coherent and every submodule of a flat $D$-module is Gorenstein flat.
\end{itemize}
\end{definition}

The notion of a \emph{Gorenstein Krull domain} (G-Krull domain for short) is also introduced in \cite{Qi-Wa17}. Examples of G-Pr\"ufer and G-Krull domains are provided in \cite{Xi22}.

Holm's metatheorem \cite{Ho04} proclaims that
\begin{itemize}
\item[] \emph{Every result in classical homological algebra has a Gorenstein counterpart.}
\end{itemize}
This sweeping statement has been largely corroborated by numerous analogies and parallels between classical and Gorenstein homological results. Yet there remain instances in which Gorenstein analogies are only partially true or require stricter hypotheses.

For example, in the classical setting, a domain $D$ is Pr\"ufer if and only if it is locally a valuation domain (\cite[Corollary~3.7.22]{Wa-Ki24a}). In the Gorenstein setting, while every G-Pr\"ufer domain is locally G-Pr\"ufer \cite[Corollary~11]{Xi23a}, the converse \emph{fails} in general \cite[Example~4.1]{Xi-Qi-Ki-Hu24}. Hence, the direct classical analogy does not hold perfectly in Gorenstein contexts.

Classically, $D$ is Pr\"ufer if and only if every overring of $D$ is flat \cite[Theorem~6.10]{La-Mc71}. Analogously, it has been shown that if $R$ is G-Pr\"ufer, then every overring of $R$ is Gorenstein flat \cite[Theorem~6]{Xi23a}. However, it remains \emph{open} whether the converse also holds. Namely, is a domain $D$ G-Pr\"ufer if every overring of $D$ is Gorenstein flat? This is a natural Gorenstein-homological parallel to a well-known Pr\"ufer criterion, but it awaits resolution.

A key feature of Dedekind domains $D$ is that \emph{every nonzero ideal is invertible}, which is equivalent to demanding that $D/I$ be quasi-Frobenius for each nonzero ideal $I$ (\cite[Theorem~13]{Xi-Zh24}). Similarly, in the Gorenstein Dedekind setting, one has:

\begin{itemize}
    \item An integral domain $D$ is G-Dedekind if and only if, for every nonzero element $a \in D$, the factor ring $D/(a)$ is  quasi-Frobenius (\cite[Theorem 11.7.7]{Wa-Ki24a}). This implies that every nonzero ideal of $D$ is a $v$-ideal. In fact, $D$ is a G-Dedekind domain if and only if it is a Noetherian divisorial domain.

    \item In general, a G-Dedekind domain is one-dimensional and Noetherian (\cite[Theorem 11.7.7 and Corollary 11.7.8]{Wa-Ki24a}), yet can exhibit different behaviors regarding \emph{integral closure}, reminiscent of how Gorenstein rings relax some classical regularity or integrality conditions.
\end{itemize}

Recent connections show that factor rings of the form $(D/xD)_w$ being quasi-Frobenius for all nonzero $x$ in $D$ characterizes G-Krull domains \cite[Theorem~2.7]{Xi24a}. Such \emph{Gorenstein-enabled} statements suggest deep analogies between multiplicative ideal theory and Gorenstein homological algebra, but also underscore certain exceptions where classical and Gorenstein theories diverge.

\begin{problem}\label{kim_prob2}
Develop explicit characterizations of Gorenstein homological properties---with an emphasis on G-Dedekind and G-Pr\"ufer domains---through both ideal-theoretic and ring-theoretic lenses. In particular:
\begin{itemize}
    \item Clarify how the $v$-operation and $\ast$-invertibility interact in the Gorenstein setting, and determine whether the $v$-operation fully captures Gorenstein analogues of classical invertibility.
    \item Investigate if and how integral closure behaves differently for G-Dedekind domains compared to classical Dedekind domains.
    \item Resolve open questions about overring flatness: under what conditions does ``every overring is Gorenstein flat'' imply that the base domain is G-Pr\"ufer?
    \item Investigate whether $\overline{D}$ is a Pr\"ufer domain when $D$ is G-Pr\"ufer.
\end{itemize}
\end{problem}

\section{Cotorsion Theory} \label{sec:cotorsion}

Cotorsion theory is a branch of homological algebra that uses the $\operatorname{Ext}$ functor to organize module categories (or abelian categories) into two complementary classes of objects, generalizing concepts such as injectivity, flatness, and their relative versions. By drawing on ideas from torsion theories and projective/injective dimensions, cotorsion theory provides a powerful framework for decomposing modules and understanding the structure of extension groups.

Throughout this section, let $R$ denote a commutative ring, possibly with zero divisors.

\begin{definition}
A \emph{cotorsion pair} $(\mathcal{A}, \mathcal{B})$ in the category of $R$-modules (or in any abelian category) consists of two classes of modules $\mathcal{A}$ and $\mathcal{B}$ such that:
\begin{enumerate}
    \item
    \[
    \mathcal{A} \;=\; \bigl\{\,M \;\bigm|\; \operatorname{Ext}^1_R(M,B)=0\text{ for all }B\in \mathcal{B}\bigr\}
    \] and
    \[
    \mathcal{B} \;=\; \bigl\{\,N \;\bigm|\; \operatorname{Ext}^1_R(A,N)=0\text{ for all }A\in \mathcal{A}\bigr\}.
    \]
    \item
    \begin{itemize}
        \item $\mathcal{A}$ is closed under direct sums, summands, and extensions.
        \item $\mathcal{B}$ is closed under direct products, summands, and extensions.
    \end{itemize}
\end{enumerate}
The classes $\mathcal{A}$ and $\mathcal{B}$ are often called \emph{orthogonal complements} of each other with respect to $\operatorname{Ext}^1_R(-,-)$.
\end{definition}

\begin{example}\label{ex:cotorsionpair}~

\begin{enumerate}
\item Let $\mathcal{M}$ be the category of all $R$-modules, $\mathscr{P}$ the class of all projective modules, and $\mathscr{I}$ the class of all injective modules. Then
\[
(\mathscr{P},\mathcal{M})
\quad\text{and}\quad
(\mathcal{M},\mathscr{I})
\]
are both cotorsion pairs in $\mathcal{M}$.

\item Over any ring $D$, the cotorsion pair $(\mathcal{A},\mathcal{B})=(\mathrm{Flat},\mathrm{Cotorsion})$ satisfies:
\[
\mathcal{A} \;=\; \{\text{flat }D\text{-modules}\},
\quad
\mathcal{B} \;=\;
\bigl\{\,N \;\bigm|\; \operatorname{Ext}^1_R(F,N)=0\text{ for all flat }F\bigr\}.
\]
\end{enumerate}
\end{example}

Cotorsion pairs underlie many useful constructions in homological algebra, including relative homological dimensions, abelian model categories, and exact categories. By focusing on how $\operatorname{Ext}^1_R(-,-)$ vanishes, one can obtain powerful decomposition theorems for modules.

\smallskip

A variant of cotorsion theory, called \emph{Tor-torsion theory}, shifts attention from $\operatorname{Ext}$ to $\operatorname{Tor}$ and hence emphasizes torsion properties and tensor products.

\begin{definition}
Let $\mathcal{A}$ and $\mathcal{B}$ be two classes of modules. The pair $(\mathcal{A}, \mathcal{B})$ is called a \textit{Tor-torsion pair} (or \textit{Tor-torsion theory}) if
\[
\mathcal{A} = \mathcal{B}^{\top}
\quad \text{and} \quad
\mathcal{B} = \mathcal{A}^{\top},
\]
where
\[
\mathcal{A}^{\top}
=
\left\{
B \in \mathcal{M}
\middle|
\operatorname{Tor}_1^R(A,B)=0 \text{ for every } A \in \mathcal{A}
\right\}.
\]
\end{definition}

In classical cotorsion theory, one requires $\operatorname{Ext}^1_R(F,M)=0$ for every flat $R$-module $F$. In Tor-torsion theory, one instead demands that $\operatorname{Tor}_1^R(F,M)=0$ for every flat $R$-module $F$. Both conditions provide different perspectives on how flat modules (or certain torsion-free classes) might interact with a given module, illuminating complementary aspects of homological dimension and torsion phenomena.

\begin{example}
Let $\mathscr{F}$ be the class of all flat $R$-modules. Then $(\mathscr{F}, \mathcal{M})$ is a (classical) Tor-torsion pair in the category of $R$-modules. In a Tor-torsion context, we focus on $R$-modules $M$ with $\operatorname{Tor}_1^R(F,M)=0$ rather than on $\operatorname{Ext}^1_R(F,M)=0$.
\end{example}

\smallskip

Cotorsion theory (and its Tor-torsion counterpart) can effectively describe structural aspects of rings that are central in multiplicative ideal theory. For instance:

\begin{enumerate}
\item
An integral domain $D$ is Dedekind if and only if the global dimension of $D$ is at most $1$. Equivalently, every nonzero ideal of $D$ is projective (i.e., invertible).

\item
An integral domain $D$ is Krull if and only if the global \emph{weak $w$-projective dimension} $w.w.gl.\dim(D)$ is at most $1$; see \cite[Theorem~4.3]{Wa-Qi19a}.
\end{enumerate}

These homological characterizations blend naturally with cotorsion-theoretic viewpoints: the condition $gl.\dim(R)\le 1$ can be recast via orthogonality requirements on $R$-modules, while the $w$-projective setting relies on extension functors in the $w$-module theory sense. Indeed, building robust (Tor-)cotorsion frameworks broadens the scope of such equivalences, reinforcing links between ring-theoretic properties (e.g., invertibility of ideals) and vanishing results in homological algebra.

\smallskip

Although there are many ways to characterize Pr\"ufer domains, cotorsion viewpoints yield elegant homological statements. First, we recall additional module-theoretic definitions:

\begin{definition}\label{def:VariousCotorsions}
Let $D$ be a domain with quotient field $K$, and let $M$ be an $D$-module.
\begin{enumerate}
\item $M$ is called a \emph{cotorsion module} if
\[
\operatorname{Ext}^1_D(F, M) \;=\; 0
\quad\text{for every flat $D$-module }F.
\]
\item $M$ is called a \emph{Warfield cotorsion module} if
\[
\operatorname{Tor}_1^D(M, F) \;=\; 0
\quad\text{for every torsion-free $D$-module }F.
\]
\item $M$ is called a \emph{Lee cotorsion module} if
\[
\operatorname{Ext}^1_D(K/D, M) \;=\; 0.
\]
\end{enumerate}
\end{definition}

\begin{theorem}[\cite{Wa-Ki24a}, Theorem~12.9.12]\label{thm:pruferCotorsion}
For a domain $D$, the following conditions are equivalent.
\begin{enumerate}
\item[(a)] $D$ is a Pr\"ufer domain.
\item[(b)] Every cotorsion $D$-module is Warfield cotorsion.
\item[(c)] Every divisible $D$-module is FP-injective.
\item[(d)] Every $h$-divisible $D$-module is FP-injective.
\item[(e)] Every Lee cotorsion $D$-module is FP-injective.
\end{enumerate}
\end{theorem}

Thus, a domain is Pr\"ufer precisely when several natural classes of modules (cotorsion, divisible, etc.) coincide homologically or admit uniform homological dimensions. This result interlinks classical ideals-based characterizations (e.g., localizing at height-one primes or checking invertibility of ideals) with module-theoretic vanishing conditions.

\smallskip

The interplay between cotorsion theory and $w$-theory highlights how star operations can refine homological statements, particularly for non-Noetherian rings or those with specialized divisorial properties. For example, focusing on $w$-flatness and $w$-cotorsion can unify certain results in multiplicative ideal theory with relative homological vanishing conditions. In this direction, the cases of P$v$MDs and Krull domains have been studied in \cite{Pu-Zh-Ta-Wa-Xi19a, Pu-Zh-Ta-Wa-Xi19b}.

\smallskip

\begin{problem}\label{kim_prob2}
Develop explicit \emph{cotorsion} and \emph{Tor-torsion} characterizations of rings and integral domains that arise in multiplicative ideal theory. In particular, investigate how homological invariants such as
\[
\mathrm{Ext}^1_R(-,-), \quad
\mathrm{Tor}_1^R(-,-), \quad
\text{global dimensions, etc.}
\]
can be used to classify such rings. Examples of interest include:
\begin{itemize}
\item {\it Pr\"ufer domains}: Extend Theorem~\ref{thm:pruferCotorsion} by identifying further cotorsion-theoretic criteria.

\item {\it Krull domains}: Relate the $w$-projective dimension viewpoint to \newline (Tor-)cotorsion approaches.

\item {\it Other star-operation domains}: Characterize specialized classes (e.g., Gorenstein domains, G-Pr\"ufer, P$v$MDs, etc.) via novel cotorsion or Tor-torsion properties.
\end{itemize}
\end{problem}

\smallskip
\section{Factorization Theory in Mori Monoids} \label{sec:factorization}
\smallskip

Factorization theory has spread into various subfields of algebra. However, in this survey we restrict ourselves to factorization theory in Mori and Krull monoids and domains. We discuss neither factorization theory in non-commutative ring theory \cite{Ba-Sm18, Sm19a, B-B-N-S23a}, nor combinatorial factorization theory with its overlap to additive combinatorics, \cite{Ge-Gr-Zh25a}, nor factorization theory in non-BF-monoids \cite{Co-Go25a}, nor factorizations into distinguished building blocks other than irreducibles \cite{Co25a}.

We start with a lemma highlighting the abundance of  Mori and Krull monoids in commutative ring theory. For more on the ideal theory of Mori monoids we refer to \cite[Chapter 5]{HK25a} and for more of this flavor in the setting of  $v$-Marot and $t$-Marot rings  see \cite{Ch-Oh22a, Ba-Po25a}.  For Mori monoids with a combinatorial flavor see \cite{F-G-R-Z24a}.
Moreover, we refer to a survey of Facchini \cite{Fa06a} for the role of Krull monoids in module theory.

\smallskip
\begin{lemma}~ \label{factorization.1}

\begin{enumerate}
\item A domain $D$ is Mori if and only if $D^{\bullet}$ is a Mori monoid if and only if $\mathcal I_v^* (D)$ is a Mori monoid.

\item If $D$ is a $v$-Marot ring, then $D$ is a Mori ring if and only if $D^{\bullet}$ is a Mori monoid.

\item Let $D$ be a domain and let $r$ be an ideal system on $D$ with $\mathcal I_r (D) \subseteq \mathcal I (D)$. Then $D$ is a Krull domain  if and only if $D^{\bullet}$ is a Krull monoid if and only if $\mathcal I_r^* (D)$ is a Krull monoid.

\item If $D$ is a Noetherian domain, then its complete integral closure $\widehat D$ is Krull and    $\mathcal I^* (D)$  is a Mori monoid.

\item Let $M$ be a Mori monoid. If $M$ is seminormal or if the  conductor $(M \colon \widehat M)$ is nonempty, then the complete integral closure $\widehat M$ is Krull. Furthermore,   the complete integral closure of a strong Mori domain is Krull.
\end{enumerate}
\end{lemma}

\begin{proof}
For (1) and (2) we refer to \cite{Ge-Ra-Re15c} and for (3) to \cite[Theorem 4.5]{Ge-Kh22b}.
The first claim in (4) is the Theorem of Mori-Nagata and the second claim follows from \cite[Example 2.1]{Ge-Ha08a}.

(5)
For a seminormal $M$ the claim follows from \cite[Theorem 5.2.8]{HK25a} and if the conductor is nonempty, then it follows from \cite[Theorem 2.3.5]{Ge-HK06a}. If $M$ is a strong Mori domain, then the claim follows from \cite[Theorems 7.7.14 and 7.10.21]{Wa-Ki24a}.
\end{proof}

We start the discussion of Arithmetic Finiteness Properties. Let $M$ be a monoid. We say that $M$ is
\begin{enumerate}
\item  {\it atomic} if every element $a \in M \setminus M^{\times}$ has a factorization into atoms.

\item {\it factorial} if one of the following two equivalent conditions hold.
      \begin{itemize}
      \item[(a)] $M$ is atomic and factorization into atoms is unique up to the order and up to associates.

      \item[(b)] $M$ is Krull and each divisorial ideal of $M$ is principal.
      \end{itemize}
\end{enumerate}
The uniqueness statement in Condition 2.(a) means that whenever $a = u_1 \ldots u_k = v_1 \ldots v_{\ell}$, with $k, \ell \in \N$ and $u_1, \ldots, u_k, v_1, \ldots, v_{\ell} \in \mathcal A (M)$, then $k=\ell$ and there is a permutation $\sigma \colon [1, k] \to [1,k]$ such that $u_i M^{\times} = v_{\sigma (i)} M^{\times}$ for all $i \in [1,k]$.

Mori monoids are atomic. Indeed, suppose that $M$ is Mori and assume to the contrary that $\Omega = \{aM \colon a \in M \setminus M^{\times} \ \text{is not a product of atoms} \}$ is nonempty. Since principal ideals are divisorial, $\Omega$ has a maximal element $aM \in \Omega$. Then there are $b , c \in M \setminus M^{\times}$ such that $a = bc$, $aM \subsetneq bM$, and $aM \subsetneq cM$. Thus, $b$ and $c$ are products of atoms and so is $a = bc$, a contradiction. The overall goal is to describe the non-uniqueness of factorizations in Mori monoids which are either not Krull or Krull with nontrivial $v$-class group. This is done by arithmetic invariants which are set in context with algebraic parameters of the respective class of Mori monoids.

We introduce length sets and some parameters controlling their structure.
If $a = u_1 \ldots u_k$, with $k \in \N$ and $u_1, \ldots,  u_k \in \mathcal A (M)$, then the number of factors $k$ is called a factorization length of $a$ and the set $\mathsf L_M (a) = \mathsf L (a)$ of all factorization lengths of $a$ is called the {\it length set} of $a$. Note that $1 \in \mathsf L (a)$ if and only if $a$ is an atom if and only if $\mathsf L (a) = \{1\}$. It is a convenient convention to set $\mathsf L (a) = \{0\}$ if $a$ is invertible. We denote by
\[
\mathcal L (M) = \{\mathsf L (a) \colon a \in M \}
\]
the {\it system of length sets} of $M$.
The monoid $M$ is said to be {\it half-factorial} if $|\mathsf L (a)|=1$ for all $a \in M$. Suppose that $M$ is not half-factorial. Then there exists some $a \in M$ with factorizations of two different lengths, say
\[
a= u_1 \ldots u_k = v_1 \ldots v_{\ell} \,,
\]
where $k, \ell \in \N$ with $k < \ell $ and atoms $u_1, \ldots, u_k, v_1, \ldots, v_{\ell}$. Then, for every $N \in \N$, we have
\[
a^N = (u_1 \ldots u_k)^i ( v_1 \ldots v_{\ell})^{N-i} \quad \text{for all $i \in [0, N]$} \,,
\]
whence $|\mathsf L (a^N)| > N$. Thus, length sets may become arbitrarily large.
For an element $a \in M$,  let $\omega_M (a) = \omega (a)$ be the smallest $N \in \N_0\cup \{\infty\}$ for which the following holds{\rm \,:}
\begin{itemize}
\item If $m \in \N$ and $v_1, \ldots, v_m \in M$ with $a \mid v_1  \ldots  v_m$, then there exists a subset $\Omega \subseteq [1,m]$  such that $|\Omega| \le N$ and $a \mid \prod_{\nu \in \Omega} v_{\nu}$.
\end{itemize}
For an atom $a \in M$, let $\mathsf t (a)$ be the smallest $N \in \N_0\cup \{\infty\}$ for which the following holds{\rm \,:}
\begin{itemize}
\item If  $m \in \N$ and $v_1, \dots , v_m \in \mathcal A(M)$ are
such that \ $a \mid v_1  \ldots  v_m$, but $a$ divides no
proper subproduct of $v_1 \ldots  v_m$, then there exist
$\ell \in \N$ and  $u_2, \ldots , u_{\ell} \in \mathcal A
(M)$ such that $v_1  \ldots  v_m = a u_2 \cdot \ldots
\cdot u_{\ell}$ and  $\max \{\ell ,\, m \} \le N$.
\end{itemize}
We say that $M$ is {\it locally tame}
if $\mathsf t (a) < \infty$ for all $a \in \mathcal A (M)$.
By definition, $\omega (a)=1$ if and only if $a$ is a prime element. If $a$ is an atom but not prime, then $\omega (a) \le \mathsf t (a)$.

The next lemma gathers some basic arithmetic properties of Mori monoids.

\smallskip
\begin{lemma} \label{factorization.2}
Let $M$ be a Mori monoid.
\begin{enumerate}
\item All length sets are finite and nonempty $($\cite[Theorem 2.2.9]{Ge-HK06a}$)$.

\item $\omega (a) < \infty$ for all $a \in M$ $($\cite[Theorem 4.2]{Ge-Ha08a}$)$.
\end{enumerate}
\end{lemma}

\smallskip
For a finite, nonempty set $L = \{m_1, \ldots, m_k\} \subset \N$ with $k \in \N$ and $m_1 < \ldots < m_k$, we call
\begin{itemize}
\item $\Delta (L) = \{ m_{i+1}-m_i \colon i \in [1, k-1] \} \subseteq \N$ the {\it set of distances} of $L$, and

\item $\rho (L) = m_k/m_1$ the {\it elasticity} of $L$.
\end{itemize}
For convenience, we set $\rho ( \{0\})=1$. Let $M$ be a Mori monoid and let $k \in \N$. Then
\begin{itemize}
\smallskip
\item $\Delta (M) = \bigcup_{L \in \mathcal L (H)} \Delta (L) \ \subseteq \N$ is the {\it set of distances} of $M$,

\smallskip
\item $\rho (M) = \sup \{ \rho (L) \colon L \in \mathcal L (H) \} \in \mathbb R_{\ge 1} \cup \{\infty\}$ is the {\it elasticity} of $M$, and

\smallskip
\item $\mathcal U_k (M) = \bigcup_{k \in L, L \in \mathcal L (M)} L \ \subseteq \N$   is the {\it union of length sets} containing $k$.
\end{itemize}
Thus, by definition,
\[
M \ \text{is half-factorial} \quad \text{if and only if} \quad \Delta (M) = \emptyset \quad \text{if and only if} \quad \rho (M)=1 \,,
\]
and, by a straightforward calculation, we get
\[
\rho (M) = \lim_{k \to \infty} \frac{\sup \mathcal U_k (M)}{k} \,.
\]
We say that $M$ satisfies the {\it Structure Theorem for Unions} if there are $k^*\geq 1$ and $N \geq 0$ such that, for all $k \ge k^*$, \ $\mathcal U_k (M)$ has the form
\[
\mathcal U_k (M) = y_k + \big( L' \cup L^* \cup L'' \big) \subset y_k + d \Z ,
\]
where  $y_k \in \Z$,  $L^*$ is a nonempty (finite or infinite) arithmetic progression with difference $d$, \ $\min L^*=0$, \ $L' \subset [-N, -1]$, and $L'' \subset \sup L^* + [1, N]$, with the convention that $L''=\emptyset$ if $L^*$ is infinite (if $M$ is half-factorial, then $\Delta (M)=\emptyset$, $d=0$, $L'=L''=\emptyset$, $L^*=\{0\}$, and $\mathcal U_k (M) = y_k+\{0\} = \{k\}$ for all $k \ge 1$).

\medskip
We say that $M$ satisfies the {\bf Main Arithmetic Finiteness Properties} if
\begin{itemize}
\item[(a)] $M$ is locally tame.

\item[(b)] The set of distances is finite.

\item[(c)] The Structure Theorem for Unions holds.

\item[(d)] The Structure Theorem for Length Sets holds (it says that all length sets are almost arithmetic multiprogressions with global bounds for all parameters; for the technical definition see \cite{Ge-HK06a}).
\end{itemize}

\medskip
Consider two simple examples of Mori domains, Cohen-Kaplansky domains and one-dimensional local Mori domains. Their arithmetic is extremely simple and it stimulates research to which extent the arithmetic of more general domains fulfills some resp. all of their arithmetic finiteness properties.
Let $D$ be a domain and let $r$ be an ideal system on $D$ such that $\mathcal I_r (D) \subseteq \mathcal I (D)$. Then $D$ is called a {\it Cohen-Kaplansky domain} if it satisfies one of the following equivalent conditions (\cite[Proposition 4.5]{Ge-Re19d})
\begin{enumerate}
\item[(a)] $D$ is atomic and has only finitely many atoms up to associates.

\item[(b)] The monoid of principal ideals is finitely generated.

\item[(c)] The monoid $\mathcal I_r^* (D)$ of $r$-invertible $r$-ideals is finitely generated.

\item[(d)] The semigroup $\mathcal I_r (D)$ of all nonzero $r$-ideals is finitely generated.
\end{enumerate}
This class of domains  was introduced by  Anderson and  Mott \cite{An-Mo92}
(see also \cite{Co-Sp12a, Cl-Go-P017a}). Every Cohen-Kaplansky domain is one-dimensional, semilocal, Noetherian with nonzero conductor.

\medskip
\begin{theorem}~ \label{factorization.3}

\begin{enumerate}
\item If $D$ is a one-dimensional local Mori domain, then $D$ satisfies the Main Arithmetic Finiteness Properties. Moreover, $M$ has finite elasticity if and only if the conductor is nonzero and the natural map $\spec ( \widehat D) \to \spec (D)$ is injective.

\item Let $M$ be a monoid such that $M_{\red}$ is finitely generated (e.g., the monoid of nonzero elements of a Cohen-Kaplansky domain). Then $M$ has finite elasticity and satisfies the Main Arithmetic Finiteness Properties.
\end{enumerate}
\end{theorem}

\begin{proof}
(1) The statement on the elasticity follows from \cite[Theorem 2.12]{An-An92} and from \cite[2.9.2+2.10.7]{Ge-HK06a}. The other statements follow from \cite[Theorem 4.1]{Ge-Go-Tr21} and from \cite[Theorem 3.9]{Ge-Ro20a}.

(2) This follows from \cite[3.1.4+4.4.11]{Ge-HK06a} and from \cite{Tr19a}.
\end{proof}

\section{Some Long-Term Arithmetic Problems} \label{sec:long-term arithmetic}

\smallskip
\subsection{Half-Factoriality}

In 1960 Carlitz proved that a ring of integers in an algebraic number field is half-factorial if and only if its class group has at most two elements. Ever since then half-factoriality is a key property studied in factorization theory and an abundance of results has been achieved (we refer to the surveys \cite{Ch-Co00, Co05a} for some early work and to \cite{Du26a} for a recent characterization of half-factoriality). We can pick only one recent highlight and one of our favorite problems.

After Carlitz's result, the question which (non-principal) orders in number fields are half-factorial received a lot of attention. Many special cases were handled (including orders in quadratic number fields and seminormal orders), but only recently Rago \cite{Ra25a} characterized half-factoriality for arbitrary orders in number fields.

In 1966 Claborn proved that every abelian group is isomorphic to the class group of a Dedekind domain. His result stimulated much further research. There is a lot of work on the realization of various types of abelian groups as class groups of Dedekind domains with wanted additional properties. We refer to some recent results and the literature cited there (\cite{Ch-Ge24a, Pe23a, Po25b}). The following problem was first posted by Steffan \cite{St86} in 1986.
It has an affirmative answer for all Warfield groups \cite{Ge-Go03} but is open in general.

\begin{problem} \label{long-term.1}
Is every abelian group isomorphic to the class group of a half-factorial Dedekind domain?
\end{problem}

\smallskip
\subsection{Finiteness of Arithmetic Invariants} A key direction of research is to characterize (in a distinguished class of monoids) the finiteness of one resp. of all the arithmetic properties as they hold true in Theorem \ref{factorization.3}. We discuss a sample of results and then we give a more precise formulation of some problems. We start with a result by Kainrath \cite{Ka05a}.

\smallskip
\begin{theorem} \label{long-term.2}
Let $D$ be a  domain that is finitely generated as a $\Z$-algebra. Then the elasticity $\rho (D)$ is finite if and only if $\mathcal C_v (\overline D)$ and $\overline D/D$ are both finite and the natural map $\spec ( \overline D) \to \spec (D)$ is injective.
\end{theorem}

We continue with finiteness results for Krull monoids which highlight the significance of their  $v$-class group for the arithmetic.

\smallskip
\begin{theorem} \label{long-term.3}
Let $M$ be a Krull monoid and let $G_0 \subseteq \mathcal C_v (M)$ denote the set of classes containing prime divisors.
\begin{enumerate}
\item If $G_0$ is finite, then the Main Arithmetic Finiteness Properties hold and $\rho (M) < \infty$.

\item Suppose that $\mathcal C_v (M)$ is infinite cyclic, say $\mathcal C_v (M)= \Z$. Then $\rho (M) < \infty$ if and only if $\Delta (M)$ is finite if and only if all $\mathcal U_k (M)$'s are finite and the Structure Theorem for Unions holds if and only if either $G_0 \cap \N$ or $G_0 \cap (-\N)$ is finite.
\end{enumerate}
\end{theorem}

\begin{proof}
The results of (1) can be found in \cite{Ge-HK06a} and (2) follows from \cite[Theorem 4.2]{Ge-Gr-Sc-Sc10}.
\end{proof}

More is known than what is given by the above results. Among others, the assumption that $G_0$ is finite can be generalized to the assumption that the Davenport constant of $G_0$ is finite.
If $\mathcal C_v (M)$ is finitely generated, then there is a characterization (in geometric terms of $G_0$) when $\rho (M)$ is finite (\cite{Gr22a}).  In particular, we do know that, in the above setting,  length sets are completely controlled by the set $G_0$ (see the next subsection), but we are far away from a characterization of when one of the Main Arithmetic Finiteness Properties (a) - (d) hold true in general.

\smallskip
\begin{problem} \label{long-term.4}
Let $M$ be a Krull monoid and let $G_0 \subseteq \mathcal C_v (M)$ denote the set of classes containing prime divisors. Characterize (in terms of $G_0$) when one resp. all of the Main Arithmetic Finiteness Properties hold.
\end{problem}

In the local case Theorem \ref{factorization.3} provides a satisfactory answer for one-dimensional domains, but only little is known for higher-dimensional local domains. Apart from the local case and apart from Krull monoids, C-monoids are the most studied class of Mori monoids with respect to the arithmetic, and they include a class of higher-dimensional Mori domains (Theorem \ref{C-monoid-3}.2). We recall their definition.

\smallskip
Let $F$ be a factorial monoid and let $M \subseteq F$ be a submonoid. Two elements $y, y' \in F$ are called $M$-equivalent (we write $y \sim y'$) if $y^{-1}M \cap F = {y'}^{-1} M \cap F$ or, in other words,
\[
\text{if for all} \ x \in F, \ \text{we have} \ xy \in M \quad \text{if and only if} \quad xy' \in M \,.
\]
$M$-equivalence defines a congruence relation on $F$ and for $y \in F$, let $[y]_M^F = [y]$ denote the congruence class of $y$. Then
\[
\mathcal C (M,F) = \{[y] \mid y \in F \} \quad \text{and} \quad \mathcal C ^* (M,F) = \{ [y] \mid y \in (F \setminus F^{\times}) \cup \{1\} \}
\]
are commutative semigroups with identity element $[1]$. $\mathcal C (M,F)$ is the {\it class semigroup} of $M$ in  $F$ and the subsemigroup $\mathcal C^* (M,F) \subseteq \mathcal C (M,F)$ is the {\it reduced class semigroup} of $M$ in $F$.

A monoid $M$ is a C-{\it monoid} if it is a submonoid of some factorial monoid $F$ such that $M^{\times} = M \cap F^{\times}$ and the reduced class semigroup is finite, and a domain $D$ is a C-domain if $D^{\bullet}$ is a C-monoid.
The next theorem gathers the main algebraic and arithmetic properties of C-monoids (see \cite{Ge-HK06a}).

\begin{theorem} \label{C-monoid-1}
Let $M$ be a {\rm C}-monoid.
\begin{enumerate}
\item $M$ is Mori, $(M \DP \widehat M) \ne \emptyset$, and $\mathcal C_v (\widehat M)$ is finite.

\item $M$ satisfies the Main Arithmetic Finiteness Properties.
\end{enumerate}
\end{theorem}

This result gives rise to the following two problems.

\begin{problem} \label{C-monoid-2a}
Let $M$ be a {\rm C}-monoid. By Theorem \ref{C-monoid-1}.(2), we know that arithmetic invariants are finite but the present methods just provide abstract finiteness results. Thus the goal is to derive the  precise values (or at least good upper bounds) for arithmetic invariants in terms of suitable algebraic invariants (say in terms of the class semigroup). Much is known is the special case when $M$ is Krull (see the next subsection). But in the general case of {\rm C}-monoids only first steps are done so far (e.g., \cite{Ge-Ka-Re15a}).
\end{problem}

\begin{problem} \label{C-monoid-2b}
Let $M$ be a Mori monoid with  $(M \DP \widehat M) \ne \emptyset$ such that $\mathcal C_v (\widehat M)$ is finite. Characterize when $M$ is a {\rm C}-monoid. The cases when $M$ is seminormal or strong Mori are of special interest (\cite{Ge-Zh19c}).
\end{problem}

There is a simple characterization of C-monoids within the class of  finitely generated monoids and a partial answer for Mori domains. But the general case is open.

\begin{theorem} \label{C-monoid-3}~

\begin{enumerate}
\item Let $M$ be a monoid such that $M_{\red}$ is finitely generated. Then $M$ is Mori with $(M \DP \widehat M) \ne \emptyset$, and  $M$ is a {\rm C}-monoid if and only if $\mathcal C_v (\widehat M)$ is finite.

\item Let $D$ be Mori domain with $(D \DP \widehat D) \ne \{0\}$ and with $\mathcal C_v (\widehat D)$ finite. If $\widehat D/ (D \DP \widehat D)$ is finite, then $D$ is a {\rm C}-domain. The converse holds if $D$ is semilocal, Noetherian but not local.
\end{enumerate}
\end{theorem}

\begin{proof}
For Part (1), we refer to \cite[Proposition 4.8]{Ge-Ha08b}. The first statement in Part (2) follows from \cite{Ge-HK06a} and the converse statement is proved by Reinhart in \cite[Corollary 4.5]{Re13a}.
\end{proof}

Every generalized Cohen-Kaplansky domain (this is an atomic domain that has only finitely many atoms, up to associates, that are not prime) is a C-domain (\cite[Proposition 4.11]{Ge-Re19d}).
Part 2 of Theorem \ref{C-monoid-3} shows that, in particular, orders in algebraic number fields are C-domains. A classical result of Furtw\"angler characterizes the ideals that can be realized as conductor ideals of an order (\cite[Theorem 2.12.10]{HK20a}). This was generalized by Reinhart to a more general setting of Noetherian domains (\cite{Re16a}), and his result gives rise to the following problem.

\begin{problem} \label{C-monoid-4}
Let $H$ be a Krull monoid with finite $v$-class group. Characterize the ideals $\mathfrak f \subseteq H$ for which there is a C-monoid $M \subseteq H$ with $\widehat M = H$ and with conductor $\mathfrak f = (M \DP \widehat M)$.
\end{problem}

C-domains have been generalized to weakly C-domains which cover a larger class of Mori domains, as given in Theorem \ref{C-monoid-3}.2. They also satisfy most of the Main Arithmetic Finiteness Properties, but it is open whether the Structure Theorem for Length Sets holds true for them.

\begin{problem} \label{C-monoid-5}
Let $M$ be a weakly C-monoid in the sense of \cite{Ka16b}. Prove or disprove that the Structure Theorem for Length Sets holds (for a first step see \cite{Ge-Ka10a}).
\end{problem}

Based on a result of Nagata, Auslander and Buchsbaum proved that every regular local Noetherian domain is factorial.
The following problem  opens the door for a wide open area where only   first steps have been done so far (e.g.,  \cite{Wi25a}).

\begin{problem} \label{long-term.5}
Study the interplay of homological properties (as discussed in the previous sections) and arithmetic finiteness properties of a domain.
\end{problem}

\smallskip
\subsection{The Transfer Krull Property} A monoid homomorphism $\theta \colon M \to B$, between atomic monoids $M$ and $B$, is called a {\it transfer homomorphism} if the following two properties hold.
\begin{enumerate}
\item[{\bf (T\,1)\,}] $B = \theta(M) B^\times$  and  $\theta^{-1} (B^\times) = M^\times$.

\item[{\bf (T\,2)\,}] If $u \in M$, \ $b,\,c \in B$  and  $\theta (u) = bc$, then there exist \ $v,\,w \in M$ such that  $u = vw$, \  $\theta (v) \in bB^{\times}$ and  $\theta (w) \in c B^{\times}$.
\end{enumerate}
It is easy to check that $\mathsf L_M (a) = \mathsf L_B \big( \theta (a) \big)$ for all $a \in M$ and $\mathcal L (M) = \mathcal L (B)$. In particular, we have  $\rho (M) = \rho (B)$ and  $\Delta (M)= \Delta (B)$.

Let $M$ be a Krull monoid  and let $G_0 \subseteq \mathcal C_v (M)$ denote the set of classes containing prime divisors. Then there is a transfer homomorphism $\theta \colon M \to \mathcal B (G_0)$, where $\mathcal B (G_0)$ is the monoid of zero-sum sequences with terms from $G_0$. This allows to study arithmetic invariants of $M$ in the combinatorial monoid $\mathcal B (G_0)$, which is studied with methods from additive combinatorics. The interplay is most strong if every class contains a prime divisor (so, if $G_0=\mathcal C_v (M)$) and $\mathcal C_v (M)$ is finite, say $\mathcal C_v (M) \cong C_{n_1} \oplus \ldots \oplus C_{n_r}$ with $1 < n_1 \mid \ldots \mid n_r$. We refer to  \cite{Sc16a,Ge-Gr-Zh25a} which demonstrates that arithmetic invariants of $M$ can be studied in terms of the group invariants $(n_1, \ldots, n_r)$.

Along the way, it was observed that there are monoids and domains, that are not Krull, but nevertheless they allow transfer homomorphisms to monoids of zero-sum sequences. Inspired by the success of the powerful methods from additive combinatorics, intense investigations have been started to figure out under which conditions we do have transfer homomorphisms to monoids of zero-sum sequences. Deep results revealed that this may happen even for non-commutative rings (\cite{Sm19a}).
Here we restrict our discussion to monoids (so, to commutative cancellative semigroups) and just mention two results. A one-dimensional local Mori domain $D$ is transfer Krull if and only if $D$ is half-factorial (whence in the trivial case only; see \cite[Theorem 5.5]{Ge-Sc-Zh17b}). Recent characterizations for orders in Dedekind domains are given in \cite{Ra25b}. For a survey and more, see \cite{Ba-Re22a, Ge-Zh20a, Ba-Po25a}.

\begin{problem} \label{long-term.6}
Take a distinguished class of Mori monoids and characterize (in terms of algebraic invariants) when they are transfer Krull.
\end{problem}

Module theory offers a rich source of not necessarily cancellative semigroups of modules that are transfer Krull (for such an example, see \cite{Ba-Sm22a}). We formulate the following problem.

\begin{problem} \label{long-term.7}
Let $R$ be a not necessarily commutative ring and let $\mathcal C$ denote a class of left $R$-modules, that is closed under direct sums, direct summands, and isomorphisms. Let $\mathcal V ( \mathcal C)$ be the semigroup of isomorphism classes of modules from $\mathcal C$, with operation induced by the direct sum. If every module from $\mathcal C$ has a semilocal endomorphism ring, then $\mathcal V (\mathcal C)$ is Krull (this was proved by Facchini; for examples, see \cite{Fa12a}). For a given class $\mathcal C$, for which $\mathcal V (\mathcal C)$ is Krull, study classes $\mathcal C'$ of modules having a transfer homomorphism $\theta \colon \mathcal V (\mathcal C') \to \mathcal V ( \mathcal C)$. For a simple (artificial) example, of how this strategy could work, see \cite[Example 5.4]{Ge-Zh20a}.
\end{problem}

\smallskip
\subsection{What can a family of  length sets look like?} In contrast to Theorem \ref{long-term.3} there is the following result by Kainrath (see \cite[Theorem 7.4.1]{Ge-HK06a}).

\begin{theorem} \label{long-term.7}
Let $M$ be a Krull monoid and let $G_0 \subseteq \mathcal C_v (M)$ denote the set of classes containing prime divisors. If $G_0$ contains an infinite abelian group, then every finite, nonempty subset of $\N_{\ge 2}$ occurs as a length set of $M$.
\end{theorem}

The same result is true for large classes of integer-valued polynomials over Krull domains and of weakly Krull Mori domains (for a sample, see \cite{Ch-Fa-Wi22a, Fa-Wi22c, Fa-Wi24a}). Surprisingly in almost all settings studied so far, monoids either satisfy the Structure Theorem for Length Sets or every finite, nonempty subset of $\N_{\ge 2}$ occurs as a length set. This dichotomy gave rise to the following problem. Let $\mathcal P_{\fin} (\N_0)$ denote the system of all finite, nonempty subsets of $\N_0$.
\begin{problem} \label{4.8}
Which families $\mathcal L \subseteq \mathcal P_{\fin} (\N_0)$ can be realized as the system of length sets of a Mori monoid?
\end{problem}

Let $\mathcal L \subseteq \mathcal P_{\fin} (\N_0)$ be a family of  finite, nonempty subsets. If $\mathcal L$ is the system of length sets of a Mori monoid, which is not a group, then the following properties hold.
\begin{itemize}
\item[(i)] $\{0\}, \{1\} \in \mathcal L $ and  all other sets of $\mathcal L $ lie in $\N_{\ge 2}$.

\item[(ii)] For every $k \in \N_0$ there is  $L \in \mathcal L $ with $k \in L$.

\item[(iii)] If $L_1, L_2 \in \mathcal L $, then there is $L \in \mathcal L $ with $L_1+L_2 \subseteq L$.
\end{itemize}
If $\mathcal L$ is a family satisfying the stronger condition that for each two $L_1, L_2 \in \mathcal L$ the sumset $L_1+L_2 \in \mathcal L$, then there is a Dedekind domain $D$ such that $\mathcal L (D) = \mathcal L$ (\cite{Ge-Zh22b}). But the general case is open.

\smallskip
The system $\mathcal L (M)$ gathers all length sets of a given monoid $M$ and Theorems \ref{long-term.3} and \ref{long-term.7} describe all the phenomena that might occur. But, in addition, one might ask after the typical structure of a length set when we pick a random element. Let $M$ be an analytic monoid. This is a Krull monoid with finite $v$-class group and with a norm function that allows to establish a theory of L-functions (\cite{Ka17a}). Then almost all (in the sense of Dirichlet density) length sets are intervals (\cite[Theorem 9.4.11]{Ge-HK06a}). Other concepts of density include the Zariski density (say for $\mathbb Q$-algebras). There is much heuristic evidence  (in a variety of settings) that {\it most} (in what sense ever) length sets are much simpler compared with the possible structure of all length sets popping up in $\mathcal L (M)$. Thus, we present the following problem.

\smallskip
\begin{problem} \label{long-term.8}
Consider a distinguished class of Mori monoids which allows us to establish a natural concept of density. Describe the structure of length sets of subsets of elements having high density.
\end{problem}

We end  with a  short discussion on semigroups of ideals that are not necessarily cancellative.  The algebraic and arithmetic structure of ideal semigroups has found wide attention in the literature (for a sample of recent contributions, see \cite{Ju12a, An-Ju-Mo19a, F-G-K-T17, Ge-Re19d,Br-Ge-Re20,Ju-Mo-Nd21,Ge-Kh22b,G-J-N-R24, Le25a}). In general, these semigroups need not be cancellative (indeed, the semigroup of all nonzero ideals of a domain is cancellative if and only if the domain is almost Dedekind \cite[Theorem 23.2]{HK98}). However, semigroups of ideals are unit-cancellative in a variety of significant cases. To discuss this situation, let $M$ be a commutative semigroup with identity. We say that $M$ is unit-cancellative if $a, u \in M$ and $au=a$ implies that $u \in H^{\times}$. Factorization theory in unit-cancellative monoids can be developed along the same lines as in the cancellative setting. In particular, atoms and length sets are defined as it is done for monoids in Section 2.

Let $D$ be a domain and let $r$ be an ideal system on $D$. Suppose that $D$ satisfies the $r$-Krull Intersection Theorem, which means that
\[
\bigcap_{n \in \N_0} (I^n)_r = \{0\} \quad \text{for all $I \in \mathcal I_r (D) \setminus \{D\}$.}
\]
Now let $I, J \in \mathcal I_r (D)$ with $(IJ)_r = I$. This implies that $(IJ^n)_r = I$ for all $n \in \N$, whence $\{0\} \ne I \subseteq \cap_{n\in \N} (J^n)_r$. Since $D$ satisfies the $r$-Krull Intersection Theorem, we infer that $J=D$, whence $\mathcal I_r (D)$ is unit-cancellative.

\smallskip
\begin{problem}
Let $D$ be a domain and let $r$ be an ideal system on $D$ such that $\mathcal I_r (D)$ is unit-cancellative. Study the system of length sets of $\mathcal I_r (D)$.
\end{problem}

\section{Appendix: Non-Commutative Rings} \label{sec:appendix}

Multiplicative ideal theory and factorization theory spread from commutative ring theory to noncommutative ring theory in the first half of the 20th century. Although this is beyond the scope of our present article, we would be remiss not even to mention this development.

A starting point of ideal theory in the non-commutative setting was the pioneering work of Asano \cite{As39a}, followed by work of Asano, Murato, and Harada \cite{As-Mu53a, Ha63a}, which led to a theory of Dedekind-like rings (including Asano prime rings and Dedekind prime rings, see \cite{Mc-Ro01a}). From the 1970s on, non-commutative Krull rings have been introduced (\cite{Ma74a, Ch81a, Je-Wa86a}), and it was Wauters who introduced non-commutative Krull monoids in order to study non-commutative semigroup algebras (\cite{Wa84a, Je-Ok07a, Ge13a}). We  refer  to some  recent articles, with a focus on surveys \cite{Ru01a, Ak-Ma16a, El19a, Ru-Ya17a,Ru25a}.

The early development of factorization theory in non-commutative rings is discussed in a survey by Smertnig \cite{Sm16a}. We briefly address two current research directions.

As it is the case in commutative ring theory, a topic of research are the questions which atomic prime rings have the finite factorization property (meaning that each element has only finitely many essentially different factorizations) resp. the bounded factorization property (meaning that all length sets are finite and nonempty; see \cite{Be-He-Le17, B-B-N-S23a}). Whereas it is quite easy to verify that a commutative Noetherian domain has bounded factorizations (see Lemma \ref{factorization.2}), it is still an open problem whether every Noetherian prime ring has bounded factorizations.

A further key problem asks which Noetherian prime rings allow a transfer homomorphism to a commutative (cancellative) monoid or even to a commutative Krull monoid.
In his seminal paper \cite{Sm13a}, Smertnig  characterized maximal orders in central simple algebras over algebraic number fields, that allow a transfer homomorphism to a commutative Krull monoid (the proof is based on theory of one-sided divisorial ideals in such orders). This result initiated much further research  (e.g., \cite{Ba-Sm15, Ba-Sm18,Sm19a,Ba22a,Ba-Sa20a}). Note that the existence of a transfer homomorphism to a commutative Krull monoid immediately implies the existence of a transfer homomorphism to a monoid of zero-sum sequences over a subset of an abelian group, which in turn opens the door for methods from additive combinatorics (see the Subsection on the Transfer Krull Property).

\medskip
\noindent
{\bf Acknowledgement.} We would like to thank our colleagues Valentin Havlovec, Giulio Peruginelli, Andreas Reinhart, and Daniel Smertnig for their feedback and all their helpful comments.

\medskip
\bibliographystyle{amsalpha}

\end{document}